\documentclass[12pt,a4paper,verbatim]{amsart}
\usepackage{amsfonts}
\usepackage{ifthen}
\usepackage{amsmath}
\usepackage{verbatim}
\usepackage{graphicx}
\usepackage{amscd,amssymb,amsthm}
\newcounter{minutes}
\divide\time by 60
\newcounter{hours}
\multiply\time by 60 \addtocounter{minutes}{-\time}

\title{Bernoulli inequality and hypergeometric functions}

\author[]{ Riku Kl\'en}
\author[]{ Vesna Manojlovi\'c }
\author[]{ Slavko Simi\'c}
\author[]{Matti Vuorinen$\dagger$}

\thanks{$\dagger$ Supported by the Academy of Finland, project 2600066611}

\address{Department of Mathematics and Statistics, University of Turku, 20014 Turku, Finland} \email{ripekl@utu.fi}
\address{Faculty of Organizational Sciences, University of Belgrade, Jove Ilica 154, 11000 Belgrade, Serbia } \email{vesnam@fon.rs }
\address{ Mathematical Institute SANU, Kneza Mihaila 36, 11000
Belgrade, Serbia} \email{ ssimic@turing.mi.sanu.ac.rs}
\address{Department of Mathematics, University of Turku, 20014 Turku,
Finland} \email{vuorinen@utu.fi}

\newtheorem{theorem}[equation]{Theorem}

\newtheorem{lemma}[equation]{Lemma}

\newtheorem{corollary}[equation]{Corollary}
\newtheorem{remark}[equation]{Remark}

\theoremstyle{remark}

\newtheorem{question}[equation]{Question}

\numberwithin{equation}{section}

\begin{document}

\begin{abstract}
Bernoulli type inequalities for functions of logarithmic type are
given. These functions include, in particular, Gaussian
hypergeometric functions in the zero-balanced case
$F(a,b;a+b;x)\,.$
\end{abstract}

\def\thefootnote{}
\footnotetext{ \texttt{\tiny File:~\jobname .tex,
          printed: \number\year-\number\month-\number\day,
          \thehours.\ifnum\theminutes<10{0}\fi\theminutes}
} \makeatletter\def\thefootnote{\@arabic\c@footnote}\makeatother

\maketitle



\bigskip

\section{Introduction}
\bigskip

The Bernoulli inequality \cite[p. 34]{mit} is often used in the
following form: For $c>1, t>0$
\begin{equation} \label{bern}
\log(1+ct) \le c \log(1+t) \,.
\end{equation}

Recently, in the study of geometric function theory, the following
variant of this classical result was proved in \cite{kmv} where it
was applied to estimate distortion under quasiconformal mappings.

\begin{theorem} \label{kmvthm} {\rm (\cite[Lemma 3.1 (7)]{kmv})}
For $0<a\le 1 \le b$ let $\varphi(t) = \max \{ t^a , t^b\}\,.$
Then for $c\ge 1$ and all $t>0$
$$
\log(1+ c \varphi(t)) \le c \max \{  \log^a(1+t), b  \log(1+t) \}
\,.
$$
\end{theorem}

Note that for $a=b=1,$ Theorem \ref{kmvthm} yields the classical
Bernoulli inequality (\ref{bern}) as a particular case.

The goal of this paper is to study various generalizations of
Theorem \ref{kmvthm}. The key problem is to find classes of
functions which are of logarithmic type so that a counterpart of
Theorem \ref{kmvthm} holds. We formulate the following question.
We write $ {\mathbb R}_+ = \{ x \in  {\mathbb R}: x >0 \}\,.$

\begin{question} \label{myq0}
For $\phi(x):=\max \{ x^a,x \} , \ 0<a<1, x\in {\mathbb R}_+$, do
there exist positive constants $c_1,c_2, c_3,c_4$ such that
\bigskip
\begin{equation} \label{myq1}
 c_1\le \displaystyle\frac{\log^p(1+\phi(x))}{ \phi(\log^p(1+x))}\le c_2, \
p>0,
\end{equation}
\begin{equation} \label{myq2}
c_3\le \displaystyle\frac{\log(1+\phi(x))\log(1+\log(1+\phi(x)))}{
\phi(\log(1+x)\log(1+\log(1+x)))}\le c_4\,?
\end{equation}
\end{question}

Our first main result is Theorem \ref{ssthm2}, which settles this
question in the affirmative.


\begin{theorem} \label{ssthm2}
The inequalities (\ref{myq1}) and  (\ref{myq2}) hold with the
constants  $c_1=(\log 2)^{p(1-a)}, c_2=1,$ $ c_3= (\log
2\log(1+\log 2))^{1-a} , c_4= 1 \,.$
\end{theorem}


The following proposition gives precise monotonicity intervals and
the proof of Theorem \ref{ssthm2} is based on it.


\begin{theorem}  \label{ssthm}
 Let $f: {\mathbb R}_+ \to {\mathbb R}_+$
 be a differentiable function and for $c\neq 0$ define
 $$
 g(x):=\displaystyle\frac{f(x^c)}{ (f(x))^c}.
 $$
We have the following
 \begin{itemize}
\item[(1)]  if $h(x):=\log(f(e^x))$ is a convex function, then
 $g(x)$ is monotone increasing for $c,x \in (0,1)$ or $c,x \in (1,\infty)$ or
 $c<0$, $x>1$ and monotone decreasing for $c \in (0,1)$, $x>1$ or $c>1, x \in (0,1)$ or $c<0$, $x \in (0,1)$;

\item[(2)] if $h(x)$ is a concave function, then $g(x)$ is
monotone increasing for $c \in (0,1)$, $x>1$ or $c>1$, $x \in (0,1)$ or $c<0$,
$x \in (0,1)$ and monotone decreasing for $c,x \in (0,1)$ or $c>1, x>1$ or
$c<0$, $x>1$.

 \item[(3)] if $h(x)$ is neither convex nor concave
on ${\mathbb R}_+$, then $g(x)$ is not monotone on ${\mathbb
R}_+$.
\end{itemize}
\end{theorem}


Next we turn our attention to the Gaussian hypergeometric
functions ${}_2F_1(a,b;c;x)\,.$ Below we also use the simpler
notation
 $F(a,b;c;x)$ omitting the subscripts. As is well-known, these functions have a logarithmic
singularity at $x=1$ for real positive triples $(a,b,c)$ with
$a+b=c\,,$ see \eqref{hypergeom}. Because of this logarithmic
behavior in the zero-balanced case $c=a+b$, it is natural to
expect that we might have a counterpart of Theorem \ref{kmvthm} in
this case, under suitable restrictions on $(a,b,c)\,.$ Our second
main result reads as follows. 

\begin{theorem} \label{2ndmain} For $c,d>0$ with $1/c+ 1/d \ge 1$ the function
defined for $r\in (0,1)$ and $p>0$ by
$$\omega(c,d,p,r) = \left( \frac{r^p}{1+r^p}\, F \left( c,d;c+d;   \frac{r^p}{1+r^p} \right) \right)^{1/p}$$
is increasing in $p\,.$ In particular,
$$
\frac{\sqrt{r}}{1+\sqrt{r}}\, F \left( c,d;c+d;
\frac{\sqrt{r}}{1+\sqrt{r}} \right) \le \left( \frac{r}{1+r}\,
F \left(c,d;c+d;   \frac{r}{1+r} \right) \right)^{1/2} \,.
$$
\end{theorem}

As we will explain in Section 3 \eqref{myobs}, this result may be
regarded as a Bernoulli type inequality for the zero-balanced
hypergeometric function. 

\section{Properties of $F(a,b;c;x)$} 

\bigskip

In this  section, we study some monotonicity properties of the
function \linebreak $F(a,b;c;x)$ and certain of its combinations
with other functions.  We first recall some well-known properties
of this function which will be used in the sequel. 

It is well known that hypergeometric functions are closely related
to the classical {\it gamma function} $\Gamma(x)$, the {\it psi
function} $\psi(x)$, and the {\it beta function} $B(x,y)$. For
Re$\thinspace x>0,~~$ Re$\thinspace y>0$, these functions are
defined by
\begin{equation} \label{myB}
\Gamma(x)\equiv \int\limits_0^{\infty}e^{-t}t^{x-1}dt, ~
\psi(x)\equiv \frac{\Gamma'(x)}{ \Gamma(x)}, ~ B(x,y)\equiv
{\frac{\Gamma(x)\Gamma(y)}{\Gamma(x+y)}},
\end{equation}
respectively (cf. \cite{as}). 
We recall the {\it difference equation} \cite[Chap. 6]{as}

\begin{equation}
\Gamma(x+1)=x\Gamma(x),
\end{equation}
and the {\it reflection property} \cite[6.1.15]{as}

\begin{equation}
\Gamma(x)\Gamma(1-x)=\displaystyle\frac{\pi}{  {\sin \pi
x}}=B(x,1-x).
\end{equation}
We shall also need the function
\begin{equation} \label{myR}
R(a,b)\! \equiv\! -2\gamma\!-\!\psi(a)\!-\!\psi(b),~ R(a)\!
\equiv\! R(a,1\!-\! a), ~ R \left( \frac{1}{2} \right) \! =\! \log
16,
\end{equation}
where $\gamma$ is the {\it Euler-Mascheroni constant} given by
\begin{equation}
\gamma=\lim\limits_{n\to \infty} \biggl( \sum_{k=1}^{n} \frac{1}{k}-
\log n   \biggr) = 0.577215\ldots.
\end{equation}

For $|x|<1$ the hypergeometric function is defined by the following series expansion
\[
  F(a,b;c;x) = \sum_{n=0}^\infty \frac{(a, n)(b,n)}{(c,n)} \frac{x^n}{n!},
\]
where $(a,0)= 1,(a,n) = \Gamma(a+n)/\Gamma(a) = a(a+1)\cdots (a+n-1), n =1,2,\dots$ is the \emph{Appell symbol} and $a,b,c \in {\mathbb R} \setminus \{0\}\,.$ 
The differentiation formula (\cite[15.2.1]{as}) reads
\begin{equation}
\frac{d}{dx}F(a,b;c;x)=\frac{ab}{c}F(a+1,b+1;c+1;x).
\end{equation}

An  important tool for our work is the following classification of
the behavior of the hypergeometric function near $x=1$ in the
three cases $ a+b < c, a+b = c,$ and $ a+b > c :$
\begin{equation}\label{hypergeom}
\left\{ \begin{array}{l} F(a,b;c;1) = \dfrac{\Gamma(c)
\Gamma(c-a-b)}{\Gamma(c-a) \Gamma(c-b)}, ~~ a+b <
c,\\[.5cm]
B(a,b)F(a,b;a\!+\!b;x)\!+\!\log(1\!-\!x) \!=\! R(a,b)\!+\!{\rm O}
 ((1\!-\!x)\log(1\!-\!x)), ~ a+b=c, \\[.4cm]
 F(a,b;c;x) = (1\!-\!x)^{c -a -b} F(c\!-\!a,c\!-\!b;c;x), ~ c < a\!+\!b.
\end{array} \right.
\end{equation}

Some basic properties of this series may be found in standard
handbooks, see for example \cite{as}. For some rational triples
$(a,b,c)$, the function $F(a,b;c;x)$ can be expressed in terms of
well-known elementary functions.
For what follows an important particular case is \cite[15.1.3]{as}
\begin{equation}\label{functiong}
  g(x) \equiv x F(1,1;2;x) = \log \frac{1}{1-x}.
\end{equation}
It is clear that for $a,b,c>0$ the function $F(a,b;c;x)$ is a
strictly increasing map from  $[0,1)$ into $[1,\infty)\,.$ For
$a,b>0$ we see by (\ref{hypergeom}) that $F(a,b;a+b;x)$ defines
an increasing homeomorphism from $[0,1)$ onto $[1,\infty)\,.$

\begin {theorem}{\cite{abrvv},\cite[Theorem 1.52]{avv}} \label{1.57}
For $a, b > 0,$ let $B=B(a,b)$ be as in $(\ref {myB})$, and let
$R=R(a,b)$ be as in $\eqref{myR}$. Then the following are true.

\begin {itemize}
\item[(1)] The function $f_1(x) \equiv \frac{ F(a,b;a+b;x)-1}{
\log (1/ (1-x))} $ is strictly increasing \smallbreak

\noindent from $(0,1)$ onto $(ab/(a+b),1/B).$

\item[(2)] The function $f_2(x) \equiv BF(a,b;a+b;x) + \log (1-x)$
is strictly decreasing from $(0,1)$ onto $(R,B).$

\item[(3)] The function $f_3(x) \equiv BF(a,b;a+b;x) + (1/x)\log
(1- x)$ is increasing from $(0,1)$ onto $(B-1,R)$ if $a,b \in
(0,1).$

\item[(4)] The function $f_3$ is decreasing from $(0,1)$ onto
$(R,B-1)$ if $a,b \in (1,\infty ).$

\item[(5)] The function
$$f_4(x) \equiv \frac{xF(a,b;a+b;x)}  {\log (1/ (1-x))}$$

is decreasing from $(0,1)$ onto $(1/B,1)$ if $a,b \in (0,1)$.

\item[(6)] If $a, b > 1,$ then $f_4$ is increasing from $(0,1)$
onto $(1,1/B).$

\item[(7)] If $a = b = 1$, then $f_4(x) = 1$ for all $x \in
(0,1).$
\end {itemize}
\end {theorem}

\bigskip

 We also need the following
refinement for some parts of Theorem \ref{1.57}.

\begin{lemma} \cite[Cor. 2.14]{pv} \label{pvlem}
For $c,d>0$, denote
$$f(x) \equiv \frac{xF(c,d;c+d;x)}  {\log (1/ (1-x))}.$$

(1) If $c\in (0,\infty)$ and $d\in (0,1/c]$, then the function $f$
is decreasing with range $(1/B(c,d),1)$;

(2) If $c\in (1/2,\infty)$ and $d\ge c/(2c-1)$, then $f$ is
increasing from $(0,1)$ to the range $(1,1/B(c,d))$.
\end{lemma}

\begin{lemma} \cite[Thm 1.5]{k} \label{kuLemma} {\it If $\max \{ a,b \} \le c$ then the coefficients of
Maclaurin power series expansion of the ratio
$F(a+1,b+1,c+1;x)/F(a,b,c;x)$ form a monotone decreasing
and convex sequence.}
\end{lemma}

%


\section{Heuristic considerations}\label{heurconsid}

\bigskip

We now apply Theorem \ref{1.57} to demonstrate that the behavior
of the hypergeometric function $F(a,b;c;x)$ in the zero-balanced
case $c=a+b$ is nearly logarithmic in the sense that some basic
identities of the logarithm yield functional inequalities for the
zero-balanced function.

Fix $x\in(0,1)$ and, for a given $p>0\,,$ a number $z \in(0,1)$
such that
$$
\log \frac{1}{1-x}= xF(1,1;2;x)= \log \left( \frac{1}{1-z}
\right)^p = p \log \frac{1}{1-z} \,.
$$
Therefore $z= 1- \sqrt[p]{1-x} \,.$

\begin{lemma} \label{mylemma2} For $c,d \in (0,1]$ define
$h(x) = x F(c,d;c+d;x)/\log(1/(1-x)), x \in[0,1)$ and let
$p\ge 1\, , B=B(c,d).$ Then for all $x\in (0,1), z= 1- \sqrt[p]{1-x},$
$$
     B\,  \ge B\, h(z) \ge B \, h(x)\,\ge 1 \quad {and} \quad F(c,d;c+d;z)\ge (1/p)F(c,d;c+d;x)\,,
$$
with equality for $c=d=1\,.$
\end{lemma}

\begin{proof} Observe that for $p\ge 1$
$$0<z=1- \sqrt[p]{1-x} \le x$$
and hence the result follows from Theorem \ref{1.57} (5). The
equality statement follows from Theorem \ref{1.57} (7).
\end{proof}

\bigskip

Next, writing the basic addition formula for the logarithm
\[
  \log z+\log w =\log(z w), \quad z,w >0 \,,
\]
in terms of the function $g$ in \eqref{functiong}, we have
\[
  g(x)+g(y) = g(x+y-xy), \quad x,y \in (0,1).
\]
Based on this observation and some computer experiments
we pose the following question:

\begin{question} \label{myq3}
(1) Fix $c,d >0$ and let  $g(x) = x F(c,d;c+d;x)$  for $x \in
(0,1)$ and set
\[
  h(x,y) = \frac{g(x)+g(y)}{g(x+y-xy)}
\]
for $x,y \in (0,1)$. For which values of $c$ and $d\,,$ is this
function bounded from below and above?

(2) Is it true that
\begin{enumerate}
  \item[a)] $h(x,y) \ge 1$, if $cd \le 1$?
  \item[b)] $h(x,y) \le 1$, if $c,d > 1$?
\end{enumerate}

(3) Can the difference
$$
d(x,y)=g(x)+g(y)-g(x+y-xy)
$$
be bounded by some constants depending on $c,d$ only?

\end{question}

\bigskip

Recall that by the Bernoulli type inequality of Theorem \ref{kmvthm} we have
\begin{equation} \label{}
\log(1+ \sqrt{r}) \le \log^{1/2}(1+ {r})
\end{equation}
for all $r \in (0,1)\,.$ In terms of \eqref{functiong} this reads
as
\begin{equation} \label{myobs}
\frac{\sqrt{r}}{1+ \sqrt{r}} \, F \left( 1,1;2;\frac{\sqrt{r}}{1+
\sqrt{r}} \right) \le \left(\frac{{r}}{1+ {r}} \, F \left( 1,1;2;\frac{{r}}{1+
{r}} \right) \right)^{1/2}
\end{equation}
for all $r \in (0,1)\,.$

\begin{question} \label{omegamonot} Fix $c,d \in (0,1]$ and let
\begin{equation} \label{}
\omega(c,d,p,r) = \left(\frac{{r}^p}{1+ {r}^p} \, F \left(
c,d;c+d;\frac{{r}^p}{1+ {r}^p} \right) \right)^{1/p}
\end{equation}
for $r\in (0,1), p >0\, . $ Is it true that for each $r\,,$
$\omega(c,d,p,x)$ is increasing in $p$? If this holds, then
\eqref{myobs} would be a special case of it.
\end{question}

The answer to this question is given in Theorem \ref{ssthm4}.

\bigskip

According to the formula \eqref{hypergeom} the function
$F(c,d;c+d;x)$ has logarithmic behavior when $x$ is close to 1.
This suggests that we may expect a Bernoulli type inequality to
hold for this function.

For what follows we fix $a,b \in (0,\infty)$ with $0 < a \le 1 \le
b$ and write for $t > 0$
\begin{equation}\label{varphi}
  \varphi(t) = \max \{ t^a,t^b \}.
\end{equation}
Next we will rewrite the inequality in Theorem \ref{kmvthm} for
the function $g$ in (\ref{functiong}) when $c=1$ and denote
\[
  g(x) = \log(1+ \varphi(r)) \equiv A
\]
implying $r = \varphi^{-1}(x/(1-x))$ and we now require, in
concert with Theorem \ref{kmvthm}, that
\[
  A \le b \max \{ \log^a (1+r),\log (1+r) \}
\]
or, equivalently,
\begin{eqnarray*}
  g(x) & \le & b \max \left\{ \log^a \left( 1+\varphi^{-1} \left( \frac{x}{1-x} \right) \right) ,\log \left( 1+\varphi^{-1} \left( \frac{x}{1-x} \right) \right) \right\}\\
  & = & b \max \left\{ g^a \left( \frac{u}{1+u} \right) , g \left( \frac{u}{1+u} \right) \right\},
\end{eqnarray*}
where $u = \varphi^{-1} (x/(1-x))$ and $g$ is given in
\eqref{functiong}. Set now $\varphi(s) = x/(1-x)$ i.e. $x =
\varphi(s)/(1+\varphi(s))$ and we have for the function $g$ in
(\ref{functiong})
\begin{equation}\label{inequalitygB}
  g \left( \frac{\varphi(s)}{1+\varphi(s)} \right) \le b \max \left\{ g^a \left( \frac{s}{1+s} \right) , g \left( \frac{s}{1+s} \right) \right\}
\end{equation}
for $s > 0$. On the basis of this discussion we ask the following
question:


\begin{question} \label{my44}
{\rm Let $c,d > 0$ and $g(x) = xF(c,d;c+d;x)$. Under which
conditions on $c$ and $d\,,$ do we have that for all $s>0$
\begin{equation}\label{inequalityg}
  g \left( \frac{\varphi(s)}{1+\varphi(s)} \right) \le \frac{b^2}{a} \max \left\{ g^a \left( \frac{s}{1+s} \right) , g \left( \frac{s}{1+s} \right) \right\},
\end{equation}
where $\varphi(s)$ is as in \eqref{varphi}? }
\end{question}

On the basis of \eqref{inequalitygB} we expect that there are
numbers $c_1,c_2 \in (0,\infty)$ such that $0 < c_1 \le 1 \le c_2$
and \eqref{inequalityg} holds for all $c,d \in (c_1,c_2)$.

\begin{question} \label{my46}
{\rm Let $g$ be as in Question \ref{my44}.
Is the following generalized version of the Bernoulli inequality true: 
\[
  g(x) \le b (1+b-a) \varphi \left( g\left( \frac{\varphi^{-1}(x/(1-x))}{1+\varphi^{-1}(x/(1-x))} \right) \right),
\]
where $\varphi(x) = \max \{ x^a,x^b \}$, $\varphi^{-1}(x) = \min
\{ x^{1/a},x^{1/b} \}$, $c,d \in (0,1)$ and $0<a<1<b$?

\bigskip

Mathematica tests show that the function
\[
  t(x) = \frac{g(x)}{b \varphi \left( g\left( \frac{\varphi^{-1}(x/(1-x))}{1+\varphi^{-1}(x/(1-x))} \right) \right)}
\]
consists of three parts: $(0,\min \{ \alpha,\beta \} )$, $( \min
\{\alpha , \beta \} , \max \{ \alpha , \beta \} )$ and $(\max \{
\alpha , \beta \} ,1)$. We easily obtain that $\alpha = 1/2$,
because then $\varphi^{-1}(x/(1-x)) = 1$. Note that $\beta$ is the solution of
  \[
    g\left( \frac{\varphi^{-1}(x/(1-x))}{1+\varphi^{-1}(x/(1-x))} \right) =1.
  \]

\begin{enumerate}
  \item[(1)] When is $\beta > 1/2$?
\end{enumerate}
 Is it true that
\begin{enumerate}
  \item[(2)] $t(x)$ is monotone on each interval $(0,\min \{ \alpha,\beta \} )$, $( \min \{\alpha , \beta \} , \max \{ \alpha , \beta \} )$ and $(\max \{ \alpha , \beta \} ,1)$?
  \item[(3)] $t(x) \ge \min \{ t(1/2),t(1-) \}$?
  \item[(4)] $t(x) \le t(\beta)$?
\end{enumerate}
}
\end{question}



\section{Answers to the questions of Section \ref{heurconsid}}


Putting $\frac{x}{1-x}=s$, we have to show that $t(s)\le b$ with
\[
t(s): = \frac{g(\frac{s}{1+s})}{ \varphi \left( g\left(
\frac{\varphi^{-1}(s)}{1+\varphi^{-1}(s)} \right) \right)}, s\in
(0,\infty).
\]

The main tool for determining the best possible bounds of $t(s)$
is given by Theorem \ref{ssthm}. Therefore we have to investigate
convexity/concavity property of the function $G(u):=\log
g \left( \frac{e^u}{1+e^u} \right)$.

\bigskip

We are in a position to formulate the following result.

\bigskip

\begin{theorem} \label{ssthm5} {  Let $c,d > 0$ and $g(x) = xF(c,d;c+d;x)\,, x \in (0,1)\,.$
 The function $G(u):=\log g \left( \frac{e^u}{1+e^u} \right)$ is concave on $(-\infty,+\infty)$
if and only if $1/c+1/d\ge 1$.}
\end{theorem}
\medskip
\begin{proof} Let us consider the function $G'$ with
$ \frac{e^u}{1+e^u}=y, \ y\in (0,1)$, and write it in
the form
$$
G'(y)=1-y+y(1-y)\frac{F'(y)}{F(y)}.
$$

Since $F'(x)=\frac{cd}{c+d} {F}(c+1,d+1;c+d+1;x)$, applying
Lemma \ref{kuLemma} we get
\begin{equation}\label{star}
\frac{F'(x)}{F(x)}=\sum_0^\infty a_nx^n,
\end{equation}
where the sequence $\{a_n\}$ is monotone decreasing and convex,
with $a_0=\frac{cd}{c+d}$.

Hence
$$
G'(y)=1+(a_0-1)y+\sum_1^\infty (a_n-a_{n-1})y^{n+1},
$$
and
$$
G''(y)=\frac{cd}{c+d}-1+\sum_1^\infty (n+1)(a_n-a_{n-1})y^n<0,
$$
since $\frac{cd}{c+d}\le 1$.

Therefore $\log g(y)$  
is concave on $(0,1)$ and, consequently, $G(u)$
is concave on $(-\infty, +\infty)$. The proof is complete.
\end{proof}

\begin{remark} \label{myrmk43}
{\rm Note that if in the proof of Theorem \ref{ssthm5},
$\frac{cd}{c+d}=1+\epsilon, \epsilon>0$ then $G''(y)$ is positive
for sufficiently small $y$ and $G$ has an inflection point since
$\lim_{y\to 1^-}G''(y)<0$. Therefore the condition $c+d\ge cd$ is
necessary and sufficient for $G$ to be concave over the
whole interval.} 
\end{remark}

\bigskip

The necessary tool for answering Questions \ref{omegamonot}-\ref{my46} is
established.

\bigskip

We shall give in the sequel a positive answer to Question
\ref{omegamonot} under the condition $1/c+1/d\ge 1$ which includes
the proposed case $c,d\in (0,1]$.

\bigskip

\begin{theorem} \label{ssthm4}  { Under the condition $1/c+1/d\ge 1$ the
function $\omega$, defined above in  Question \ref{omegamonot}, is
monotone increasing in $p$.}
\end{theorem}

\begin{proof}
Denote equivalently
$$
\omega(p)=\left( g \left( \frac{e^{pt}}{1+e^{pt}} \right) \right)
^{1/p}, t<0, p>0,
$$
where $g(x):=x F(c,d,c+d;x)=xF(x)$.

We get
$$
\frac{\omega'}{\omega}=(\log\omega)'= \left(
\frac{\log(g(\frac{e^{pt}}{1+e^{pt}}))}{p}
\right)'=\frac{\Omega(p)}{p^2},
$$
with
$$
\Omega(p):=pt\frac{e^{pt}}{(1+e^{pt})^2}\frac{g'\left(
\frac{e^{pt}}{1+e^{pt}} \right)}{g\left( \frac{e^{pt}}{1+e^{pt}}
\right)}-\log g\left( \frac{e^{pt}}{1+e^{pt}} \right).
$$

Changing variable $\frac{e^{pt}}{1+e^{pt}}:=x, x\in(0,1/2)$ and
recalling the definition of $g$, we obtain
\begin{eqnarray*}
\Omega(x) & = & x(1-x)\left(\frac{F(x)+xF'(x)}{xF(x)}\right)\log\frac{x}{1-x}-\log(xF(x))\\
& = & \left( 1-x+x(1-x)\frac{F'(x)}{F(x)} \right) \log\frac{x}{1-x}-\log x-\log
F(x).
\end{eqnarray*}

 From the proof of Lemma \ref{kuLemma}, we derive the following
inequalities
$$
(1-x)\frac{F'(x)}{F(x)}<a_0; \ \ \log F(x)<-a_0\log(1-x).
$$

Noting that $\log\frac{x}{1-x}<0$ for $x\in(0,1/2)$, we get

\begin{eqnarray*}
\Omega(x) & > & (1-x+a_0 x)\log\frac{x}{1-x}-\log x+a_0\log(1-x)\\
& = & (a_0-1)\left(x\log\frac{x}{1-x}+\log(1-x)\right)\\
& = &
(1-a_0)\left(x\log\frac{1}{x}+(1-x)\log\frac{1}{1-x}\right)\ge 0,
\end{eqnarray*}
since
$$
1-a_0=1-\frac{cd}{c+d}=\frac{cd}{c+d}\left(\frac{1}{c}+\frac{1}{d}-1\right)\ge
0.
$$

Therefore $\omega'>0$ and $\omega(p)$ is monotone increasing, as
required.

\end{proof}

\begin{remark} \label{finrmk1} {\rm It is evident from the proof of Theorem \ref{ssthm4}
that the function $\omega(p):=(g(x^p/(1+x^p)))^{1/p}, 0<x<1, p>0$
is monotone increasing in $p$ for any $g(x)=x_2F_1(a,b,c;x),
a,b,c>0$ and $ab\le c, \max \{ a,b \} \le c$.

The same is valid for the conclusion of Theorem \ref{ssthm5}.

Note also that the function $\omega(p)$ is not monotone in $p$ for
$x>1$. For example in the case $(c,d,x) = (1,1,4)$, when $g$ is as
in \eqref{functiong}, we  obtain $\omega(1) = \log 5 \approx
1.61$, $\omega(2) = (\log 17)^{1/2} \approx 1.68$ and $\omega(4) =
(\log 257)^{1/4} \approx 1.53$.}
\end{remark}

\bigskip

\begin{remark} {\rm From Theorem \ref{ssthm5} it follows that the
function $g(e^x/(1+e^x))$ is log-concave on $\mathbb R$. In
particular, the function $g(x^p/(1+x^p))$ is log-concave in $p$,
that is
$$
g\left(\frac{x^p}{1+x^p}\right)g\left(\frac{x^q}{1+x^q}\right)\le
g^2\left(\frac{x^{(p+q)/2}}{1+x^{(p+q)/2}}\right), x>0,
p,q\in\mathbb R.
$$

\bigskip

Also, from Theorem \ref{ssthm4} we get that the function
$\frac{\log g\left(\frac{x^p}{1+x^p}\right)}{p}$ is monotone
increasing in $p$, that is $\log g(\frac{x^p}{1+x^p})$ is
sub-additive on $\mathbb R_+$.

Hence,
$$
g\left(\frac{x^p}{1+x^p}\right)g\left(\frac{x^q}{1+x^q}\right)\le
g\left(\frac{x^{p+q}}{1+x^{p+q}}\right), p,q>0, 0<x<1.
$$

\bigskip

Both corollaries are valid for the class of functions $g$ defined
in Remark \ref{finrmk1}. }
\end{remark}

\bigskip

An answer to part (1) of Question \ref{my46} is given by the following
assertion.

\begin{theorem} \label{ssthm55}  Let $\beta$ be as in Question \ref{my46}. We have that

\begin{enumerate}
  \item[(1)] $\beta > 1/2$ if $(a_0-1)/h\le c_0$;
  \bigskip
  \item[(2)] $\beta<1/2$ if $(a_0-1)/h\ge c_1$,
\end{enumerate}
where
$$
c_0=1-\frac{1}{2\log 2}\approx 0.27865; \ c_1=\frac{1}{\log
2}-1\approx 0.4427\,\,
{\text
and}\,\,
a_0=\frac{cd}{c+d}, \ h=\frac{a_0^2}{c+d+1}.
$$
\end{theorem}
\bigskip
\begin{proof} (1) Firstly note that the functions $g, \varphi$ and $\varphi^{-1}$  are strictly
increasing. Therefore $\beta>1/2$ if $g(1/2)<1$ and {\rm vice
versa}.

From the relation \eqref{star} we obtain

$$
(1-t)\frac{F'(t)}{F(t)}= a_0+\sum_1^\infty (a_n-a_{n-1})t^n\le
a_0-(a_0-a_1)t,
$$
since $\{a_n\}$ is a monotone decreasing sequence.

Therefore,
$$
\frac{F'(t)}{F(t)}\le a_0 \frac{1}{1-t}-(a_0-a_1)\frac{t}{1-t},
$$
and, integrating over $[0,x]$, we get
$$
\log F(x)-\log F(0)\le -a_0\log(1-x)+(a_0-a_1)(x+\log(1-x)).
$$

Putting $x=1/2$, one can see that the condition $g(1/2)\le 1$ is
satisfied if
$$
a_0\log 2+(a_0-a_1)(1/2-\log 2)\le \log 2
$$
i.e.,
$$
\frac{a_0-1}{a_0-a_1}\le 1-\frac{1}{2\log 2}=c_0.
$$

By the above remark we have that in this case $\beta>1/2$.


(2) Since $\{a_n\}$ is a convex sequence we conclude that
$\{a_{n-1}-a_n\}$  is monotone decreasing sequence.

Hence
\begin{eqnarray}
  (1-t)\frac{F'(t)}{F(t)} & = & a_0-\sum_1^\infty (a_{n-1}-a_n)t^n\label{eq*}\\
  & \ge & a_0-(a_0-a_1)t(1+t+t^2+\cdots) =a_0-(a_0-a_1)\frac{t}{1-t}.\nonumber
\end{eqnarray}

Therefore,
$$
\frac{F'}{F}\ge
(2a_0-a_1)\frac{1}{1-t}-(a_0-a_1)\frac{1}{(1-t)^2},
$$
and, integrating over $t\in [0, x]$, we get
$$
\log F(x)\ge -(2a_0-a_1)\log(1-x)-(a_0-a_1)\frac{x}{1-x}.
$$

\bigskip

Putting there $x=1/2$ we see that the condition $g(1/2)\ge 1$ is
satisfied if
$$
(2a_0-a_1)\log 2-(a_0-a_1)\ge \log 2,
$$
which is equivalent with
$$
\frac{a_0-1}{a_0-a_1}\ge \frac{1}{\log 2}-1=c_1.
$$

\bigskip

From \eqref{eq*} we get
$$
(1-t)F'(t)=F(t) \left( a_0-\sum_1^\infty (a_{n-1}-a_n)t^n \right) ,
$$
i.e.,
$$
\frac{cd}{c+d} F(c+1,d+1,c+d+1;t)=F(c,d,c+d;t) \left(
a_0-\sum_1^\infty (a_{n-1}-a_n)t^n \right) ,
$$
and, comparing power series coefficients, we easily obtain
$$
a_0=\frac{cd}{c+d}; \ \ a_0-a_1=h=\frac{c^2 d^2}{(c+d)^2(c+d+1)}.
$$
\end{proof}

\begin{corollary} \label{my49} We have that $g(1/2)<1$ if $a_0\le 1$.
\end{corollary}

\begin{remark}\label{gamma} {\rm Note that the condition $cd\le 1$ implies
$1/c+1/d\ge 1$, that is $a_0\le 1$. Therefore Theorems \ref{ssthm4} and \ref{ssthm}
could be applied to the expression $T(s)$. Moreover, by Corollary
\ref{my49} we have
$$
g(1/2)<1=g \left( \frac{ \gamma}{1+\gamma} \right),
$$
where $\gamma$ is the unique solution of the equation
$g(\frac{s}{1+s})=1$, and, since $g$ is a monotone increasing function, we conclude that
$\gamma>1$.}
\end{remark}

The following assertion is a
counterpart of Theorem \ref{ssthm4}.

\bigskip

\begin{lemma}\label{logconlemma}
For fixed $s>0$ and $c,d>0$ with $cd\le 1$, the function $\frac{g(\frac{s^p}{1+s^p})}{p}$ is
monotone decreasing in $p, \ p\in(0,\infty)$.
\end{lemma}

\begin{proof} \ Denote equivalently
$$
w(p)= \frac{g(\frac{e^{pt}}{1+e^{pt}})}{p}, \ p>0; \ t\in\mathbb
R.
$$
We have
$$
p^2w'(p)=pt\frac{e^{pt}}{(1+e^{pt})^2}g'(\frac{e^{pt}}{1+e^{pt}})-g(\frac{e^{pt}}{1+e^{pt}}):=A(p)
$$

Changing variable $\frac{e^{pt}}{1+e^{pt}}=x$, we get

$$
A(x)=x(1-x)\log{\frac{x}{1-x}} g'(x)-g(x), \ 0<x<1.
$$

Now, Lemma \ref{pvlem} tells us that the function
$\frac{g(x)}{-\log(1-x)}$ is monotone decreasing if $cd\le 1$ and
this is equivalent to $g(x)\ge -(1-x)\log(1-x) g'(x)$.

Therefore,
$$
A(x)\le
(1-x)g'(x)(x\log\frac{x}{1-x}+\log(1-x))=-(1-x)g'(x)(x\log\frac{1}{x}+(1-x)\log\frac{1}{1-x})\le
0
$$
since $g$ is an increasing function.

\end{proof}

An immediate consequence is the next corollary.

\bigskip

\begin{corollary}  \label{logcor} For $b\ge a>0$ and $cd\le 1$, the inequality
$$
1\le \frac{g(\frac{s^b}{1+s^b})}{g(\frac{s^a}{1+s^a})}\le
\frac{b}{a}
$$
holds for arbitrary $s\ge 1$.
\end{corollary}

\bigskip

Another interesting result follows from Lemma \ref{logconlemma}.

\bigskip

\begin{corollary} \label{logcor1} \ For $c,d>0$ with $cd\le 1$,
the following inequality
$$
g(\frac{s^p}{1+s^p})+g(\frac{s^q}{1+s^q})\ge
g(\frac{s^{p+q}}{1+s^{p+q}})
$$
holds true for arbitrary $s,p,q>0$.
\end{corollary}

%
%

An answer to Question \ref{my44} with improved constant is given in
the next theorem.

\begin{theorem}
Let $0<a\le 1\le b < \infty$ and let $\varphi(s)$ be defined as in \eqref{varphi} and $c,d>0$ with $cd\le 1$.
Then the inequality

\begin{equation}
  g \left( \frac{\varphi(s)}{1+\varphi(s)} \right) \le \frac{b}{a} \max \left\{ g^a \left( \frac{s}{1+s} \right) , g \left( \frac{s}{1+s} \right) \right\},
\end{equation}
holds for each $s>0$.
\end{theorem}

\begin{proof} Analyzing the structure of the above inequality, we
decide that our task is to find an upper bound for the expression
$T(s)$ given by
$$
T(s)=\begin{cases}
g \left( \frac{s^a}{1+s^a} \right)\left/ g^{a}\left( \frac{s}{1+s}\right) \right., & 0<s\le 1;\\
g\left( \frac{ s^b}{1+s^b}\right)\left/ g^{a}\left( \frac{s}{1+s}\right)\right., & 1<s\le \gamma;\\
g\left( \frac{ s^b}{1+s^b}\right) \left/ g\left( \frac{s}{1+s}\right)\right., & \gamma<s,
\end{cases}
$$
where $\gamma$ is the unique solution of the equation
$g(\frac{s}{1+s})=1$. By Remark \ref{gamma}, $\gamma > 1$.

Applying the second part of Theorem \ref{ssthm} we see that $T(s)$ is
monotone decreasing for $s\in(0,1)$. Therefore, in this case we
have
$$
T(s)<\lim_{s\to 0^+} T(s)=\lim_{s\to
0^+}\frac{s^a/(1+s^a)}{(s/(1+s))^a}\frac{F(c,d,c+d;
s^a/(1+s^a))}{(F(c,d,c+d; s/(1+s)))^a}=1.
$$

Now, for $1<s\le \gamma$, write
$$
T(s)=
\frac{g(\frac{s^a}{1+s^a})}{g^a(\frac{s}{1+s})}\frac{g(\frac{s^b}{1+s^b})}{g(\frac{s^a}{1+s^a})}=T_1(s)T_2(s).
$$

By Theorem \ref{ssthm}, $T_1(s)=
\frac{g(\frac{s^a}{1+s^a})}{g^a(\frac{s}{1+s})}$ is monotone
increasing in $s$.

Therefore, for $1<s\le \gamma$, we get
$$
T_1(s)\le T_1(\gamma)= g(\frac{\gamma^a}{1+\gamma^a})\le
g(\frac{\gamma}{1+\gamma})=1,
$$
and
\begin{equation}\label{T2ba}
T_2(s)= \frac{g(\frac{s^b}{1+s^b})}{g(\frac{s^a}{1+s^a})}\le
\frac{b}{a}
\end{equation}
by Corollary \ref{logcor}.

\bigskip

Analogously, in the case $s>\gamma$ we have
\begin{equation}\label{Tb}
T(s)= \frac{g(\frac{s^b}{1+s^b})}{g(\frac{s}{1+s})}\le b.
\end{equation}

\bigskip

The assertion follows from \eqref{T2ba} and \eqref{Tb}.

\end{proof}

Finally, an answer to Question \ref{my46} is given in the next

\begin{theorem} \label{ssthm7} \ The inequality
\[
  g(x) \le b \
  \varphi \left( g\left( \frac{\varphi^{-1}(x/(1-x))}{1+\varphi^{-1}(x/(1-x))} \right) \right),
\]
holds for $x\in (0,1)$, where $\varphi(y) = \max \{ y^a,y^b \}$,
$\varphi^{-1}(y) = \min \{ y^{1/a},y^{1/b} \}$, $0<a<1<b$ and
$c,d>0, \ cd\le 1.$
\end{theorem}

\begin{proof}
Changing variable $\frac{x}{1-x}=\varphi(s), \ s\in \mathbb R_+$,
we get
$$
t(s)=\begin{cases}
g \left( \frac{s^a}{1+s^a} \right) \left/ g^{a}\left( \frac{s}{1+s}\right)\right., & 0<s\le 1;\\
g\left( \frac{s^b}{1+s^b}\right)  \left/ g^{a}\left( \frac{s}{1+s}\right)\right., & 1<s\le \gamma;\\
g\left( \frac{s^b}{1+s^b}\right) \left/ g^{b}\left( \frac{s}{1+s}\right) \right., & \gamma<s,
\end{cases}
$$
where $\gamma$ is the unique solution of the equation
$g(\frac{s}{1+s})=1$.

\bigskip

Proceeding similarly as above, we obtain
$$
t(s)<\lim_{s\to 0^+}t(s)=1,
$$
in the case $0<s\le 1$;

\bigskip

for $1<s\le \gamma$, we have

$$
t(s)= \frac{g \left( \frac{s^b}{1+s^b} \right) }{g \left( \frac{s}{1+s} \right) }
g^{1-a} \left( \frac{s}{1+s} \right) \le b \ g^{1-a} \left( \frac{\gamma}{1+\gamma} \right) =b,
$$
by Corollary \ref{logcor};

\bigskip

 for $s>\gamma$, by Theorem \ref{ssthm} we get
$$
t(s)<t(\gamma)=
\frac{g \left( \frac{\gamma^b}{1+\gamma^b} \right) }{g \left( \frac{\gamma}{1+\gamma} \right) }\le
b,
$$
by Corollary \ref{logcor} again.

\bigskip

Therefore $t(s)\le b$ which proves Theorem \ref{ssthm7}.

\end{proof}

\section{Proofs of the main theorems}
In this section we give the proofs for the main theorems.

\begin{proof}[{\bf Proof of Theorem \ref{ssthm}.}]

We shall prove part (1) only. The proof of part (2) is similar and
the assertion of (3) follows from the former considerations.

Since $h$ is convex, $h'=\frac{e^x f'(e^x)} {f(e^x)}$ is an
increasing function, that is, if $u>v$ then

\begin{equation} \label{mystar}
 \displaystyle\frac{ e^u f'(e^u)}{ f(e^u)} >
\displaystyle\frac{e^v f'(e^v)} { f(e^v)}.
\end{equation}
Now,
\begin{eqnarray*}
\displaystyle\frac{g'}{ g} & = & c \left(
\displaystyle\frac{x^{c-1}f'(x^c)}{ f(x^c)}-
                    \displaystyle\frac{f'(x)}{f(x)} \right)\\
& = & ce^{-\log x} \left( \displaystyle\frac{e^{c\log x}
f'(e^{c\log x})}{ f(e^{c\log x})}- \displaystyle\frac{e^{\log x}
f'(e^{\log x})}{ f(e^{\log x})}
             \right),
\end{eqnarray*}
and by $(\ref{mystar})$, the conclusion of the part (1) follows by
comparing $c\log x$ with $\log x$.
\end{proof}

\bigskip

Applying Theorem \ref{ssthm} we are able to give an answer to the
above Question \ref{myq0} and prove Theorem \ref{ssthm2}. Before that, we introduce the following lemma.

\begin{lemma} \label{mylemma1}
\begin{enumerate}
\item The expression $w(x):=e^x+r(r(x))(e^x-1-r(x)),  r(x)=
\log(1+e^x),$ is positive for $x\in \mathbb R$. \item The function
$v(x):=r(r(x))=\log(1+\log(1+e^x))$ is log-concave.
\end{enumerate}
\end{lemma}

\begin{proof} (1) Putting $1+e^x=e^t, t>0$, we obtain
\begin{eqnarray*}
  w(t) & = & e^t-1+\log(1+t)(e^t-2-t)=(e^t-1)(1+\log(1+t))-(1+t)]\log(1+t)\\
  & > & t(1+\log(1+t))-(1+t)\log(1+t)=t-\log(1+t)>0.
\end{eqnarray*}

(2) By differentiation, we get
$$
\frac{d^2 \log(v(x))}{d x^2} = \frac{vv''-(v')^2}{v^2}=-
\displaystyle\frac{e^x w(x)}{ ((1+e^x)(1+\log(1+e^x)))^2},
$$
which is negative by (1).
\end{proof}


\begin{proof}[Proof of Theorem \ref{ssthm2}.]
Since $\phi(u)=u^a$ for $0<u\le 1$ and $\phi(u)=u$ if $u\ge 1$, we
easily get

\begin{equation}
 f_1(x)=\displaystyle\frac{\log^p(1+\phi(x))}{ \phi(\log^p(1+x))}=
 \begin{cases}
 \frac{\log^p(1+x^a)}{ (\log^p(1+x))^a} , & 0<x\le 1;  \\
               (\log(1+x))^{p(1-a)},  & 1<x\le e-1;\\
                1, & x > e-1.
\end{cases}
\end{equation}

For the proof of \eqref{myq1}, we will apply Theorem \ref{ssthm}
and  show first that the function $r(x):=\log(1+e^x)$ is
log-concave on $\mathbb R$.
\bigskip
Indeed, since
$$
r'(x)=\displaystyle\frac{e^x}{ 1+e^x}, \ \
r''(x)=\displaystyle\frac{e^x}{ (1+e^x)^2},
$$
we get
$$
rr''-(r')^2= \displaystyle\frac{e^x}{
(1+e^x)^2}(\log(1+e^x)-e^x)<0,
$$
because $\log(1+t)<t, \ t>0$.
\bigskip
Since $r(x)$ is log-concave, Theorem \ref{ssthm} gives
$$
f_1(1)=(\log 2)^{p(1-a)}\le f_1(x)<1=\lim_{x\to 0} f_1(x), \ x\in
(0,1].
$$
Also,
$$
f_1(1)=(\log 2)^{p(1-a)}< f_1(x)\le 1= f_1(e-1), \ x\in (1,e-1].
$$
\bigskip
Hence, $c_1=(\log 2)^{p(1-a)}, c_2=1$ and those bounds are best
possible.

Answering \eqref{myq2}, we proceed analogously. Denote
$$s(x):=\log(1+x)\log(1+\log(1+x))$$
and let $x_0, \ x_0\approx 2.4555$ be the unique positive solution
of the equation $s(x)=1$.
\bigskip
By the definition of $\phi$, we get
\begin{equation}
f_2(x)=\begin{cases}  \frac{ s(x^a)}{ (s(x))^a}, & 0<x\le 1;\\
               (s(x))^{1-a}, & 1<x\le x_0;\\
                1, & x > x_0.
        \end{cases}
\end{equation}

By Lemma \ref{mylemma1} (2)
 it is evident that $s(e^x)=r(x)r(r(x))=r(x)v(x)$ is a log-concave function
  since it is represented by the
product of two log-concave functions.

 Applying the second part of
Theorem \ref{ssthm}, we get
$$
\displaystyle\frac{s(x^a)}{ (s(x))^a} < \lim_{x\to
0}\displaystyle\frac{s(x^a)}{ (s(x))^a}=1;
$$
$$
\displaystyle\frac{s(x^a)}{ (s(x))^a}\ge \displaystyle\frac{s(1)}{
(s(1))^a}=(\log 2\log(1+\log 2))^{1-a},
$$
for $0<x\le 1$.
\bigskip
Since $s(x)$ is an increasing function, it follows that
$$
(s(1))^{1-a}<f_2(x)\le (s(x_0))^{1-a}=1,
$$
for $1<x\le x_0$.
\bigskip
Hence for $x>0$,
$$
(\log 2\log(1+\log 2))^{1-a}\le f_2(x)\le 1,
$$
and those bounds are best possible.
\end{proof}

\begin{remark} {\rm Although Question \ref{myq0} can be solved by the method of
\cite[Lemma 3.1]{kmv}, an application of Theorem \ref{ssthm} gives
the result more efficiently. }
\end{remark}

\medskip

\begin{remark} {\rm An affirmative answer to Question \ref{myq3} is given
in \cite{sv}.}
\end{remark}





\end{document}